\documentclass[a4paper,12pt]{amsart}

\usepackage{amsmath,amssymb,amsthm,latexsym,mathrsfs,dsfont}

\newtheorem{theorem}{Theorem}
\newtheorem{corollary}[theorem]{Corollary}
\newtheorem{proposition}[theorem]{Proposition}
\newtheorem{definition}[theorem]{Definition}

\newcommand{\m}{\mathcal} 

\begin{document}

\title[Subspaces Discerning Nullcontinuity]{Subspaces Discerning Nullcontinuity}
         
\author[M.\ Thill]{Marco Thill}

\address{BD G.\,-\,D.\ Charlotte 53,
                 L -\,1331 Luxembourg\,-\,City,
                 Europe}

\email{math@pt.lu}

\date{December 10, 2004}

\subjclass[2000]{28C05}

\keywords{nullcontinuous, null-continuous, subspace}

\begin{abstract}
Given positive linear functional $\ell$ on a vector lattice $\mathscr{L}$ of
real functions, and a vector subspace $\mathcal{M}$ of $\mathscr{L}$, we
construct a vector subspace $\mathcal{P}(\mathcal{M})$ of $\mathcal{M}$ in
such a way that 1) $\ell$ is nullcontinuous on $\mathcal{P}(\mathcal{M})$, and
2) if $\ell$ is nullcontinuous on $\mathcal{M}$ then $\mathcal{P}(\mathcal{M})$
is all of $\mathcal{M}$. We mention here that this result continues to hold for
quite general modes of convergence, including $\tau$\,-\,continuity. Our
construction uses a new method involving the ``kernel'' of a seminorm.
\end{abstract}

\thanks{My thanks go to Torben Maack Bisgaard for critical reading of the manuscript.}

\maketitle


\section{Basic notation and terminology}

We consider a \underline{fixed} structure $( X, \mathscr{L}, \ell, \rho )$
consisting of a non-empty set $X$, a vector sublattice $\mathscr{L}$
of $\mathds{R}^X$, a positive linear functional $\ell$ on $\mathscr{L}$,
and a set $\rho$ of non-empty downward directed subsets $I$ of
$\mathscr{L}_+$ with $\bigwedge I = 0$ in $\mathds{R}^X$.

\smallskip
The symbols $\m{M}$, $\m{M}_1$, $\m{M}_2$ shall henceforth
denote \underline{variable} vector subspaces of $\mathscr{L}$.

\begin{definition}\label{rhocont}
We shall say that $\ell$ is \underline{$\rho$\,-\,continuous on $\m{M}$}, if
\[ \bigwedge _{f \in I} \,\ell (f) = 0 
\textit{ for all } I \in \rho \textit{ with } I \subset \m{M}. \]
\end{definition}

If $\rho = \{ \,\{ \,f_n \,\} \subset \mathscr{L}_+ : 
f_n \downarrow 0 \text{ in } \mathds{R}^X \,\}$ 
then $\rho$\,-\,continuity is the same as nullcontinuity. 
If $\rho$ is the set of \emph{all} non-empty downward directed
subsets $I$ of $\mathscr{L}_+$ with $\bigwedge I = 0$ in $\mathds{R}^X$,
then $\rho$\,-\,continuity is also known as $\tau$\,-\,continuity.

\section{Review of terminology concerning vector lattices}

The vector space $\m{M}$ is called \underline{solid} in $\mathscr{L}$,
if the conditions $f \in \mathscr{L}$, $g \in \m{M}$, and $| f | \leq | g |$
together imply that $f$ belongs to $\m{M}$, cf.\ e.g.\ \cite[1.3.9]{CFW}.
In this case $\m{M}$ is a vector sublattice of $\mathscr{L}$,
because then $| f | \in \m{M}$ whenever $f \in \m{M}$.

\smallskip
A seminorm $q$ on $\mathscr{L}$ is called a \underline{lattice seminorm},
if for all $f,g \in \mathscr{L}$ with $| f | \leq | g |$ one has $q(f) \leq q(g)$,
cf.\ \cite[1.10.1]{CFW}. (An equivalent requirement is that $q(f) = q ( | f | )$
for all $f \in \mathscr{L}$.) In the affirmative case, the ``kernel''
$\{\,f \in \mathscr{L} : q(f) = 0\,\}$ of the lattice seminorm $q$ is a solid
subspace of $\mathscr{L}$, and thereby a solid vector sublattice of $\mathscr{L}$.

\smallskip
The vector lattice $\mathscr{L}$ is called \underline{Stonean}, if it contains with
each function $f$ the function $f \wedge 1_X$, cf.\ e.g.\ \cite[2.5.14]{CFW}.
This property is inherited to all solid vector subspaces.

\section{The notion of $\rho$\,-\,regularity: \\
definition and simple properties}

\begin{definition}
Denote by $T$ the set of functions $g : X \to [0,1]$, such that
$gf \in \mathscr{L}_+$ for all $f \in \mathscr{L}_+$. (The plus signs are a
matter of convenience.)
\end{definition}

One checks that $T$ is a convex sublattice of $[0,1]^X$.
Clearly $0$, $1_X \in T$. If $g \in T$, so is $1_X-g$, and both
$g$ and $1_X - g$ are positive, so multiplication by them
maps $\mathscr{L}_+$ to itself in an order preserving way.

\begin{definition}\label{S}
We put
\[ S(\m{M}) := \{ \,g \in T : \bigwedge _{f \in I} \,\ell (gf) = 0 
\textit{ for all } I \in \rho \textit{ with } I \subset \m{M} \,\}. \]
This is a convex sublattice of $[0,1]^X$. If $\m{M}_1 \subset \m{M}_2$,
then $S(\m{M}_2) \subset S(\m{M}_1)$.
\end{definition}

\begin{proof}
We note first that the ``$\bigwedge$'' in the defining
relation for $S(\m{M})$ actually is a limit since each $I \in \rho$ is
downward directed. This implies that $S(\m{M})$ is convex and
thence a sublattice of $[0,1]^X$.
\end{proof}

\begin{definition}\label{defreg}
We shall say that $\ell$ is \underline{$\rho$\,-\,regular
on $\m{M}$}, if for every $h \in \m{M}$ one has
\[  \ell ( | h | ) = \bigvee _{g \in S(\m{M})} \,\ell ( g | h | ). \]
\end{definition}

(This can be put in terms of $h \in \m{M}_+$ in
case $\m{M}$ is a vector sublattice of $\mathscr{L}$.)
We shall reformulate this condition in the next two items.

\begin{definition}\label{q}
A lattice seminorm $q_{\m{M}}$ is defined on $\mathscr{L}$ by putting
\[ q_{\m{M}}(h): = \bigwedge _{g \in S(\m{M})} \,\ell ( ( 1_X - g ) | h | ) \] 
for every $h \in \mathscr{L}$. If $\m{M}_1 \subset \m{M}_2$,
then $q_{\m{M}_1} \leq q_{\m{M}_2}$.
\end{definition}

\begin{proof}
The ``$\bigwedge$'' is a limit as $S(\m{M})$ is upward
directed, so $q_{\m{M}}$ is a semi\-norm. The last part
follows from the last part of definition \ref{S}.\end{proof}

\begin{theorem}\label{regular}
The functional $\ell$ is $\rho$\,-\,regular on $\m{M}$
if and only if $q_{\m{M}}$ vanishes identically on $\m{M}$.
\end{theorem}

\begin{proof} Let $h \in \m{M}$.
Using that $S(\m{M})$ is upward directed, one finds
\begin{align*}
q_{\m{M}} (h) & = \bigwedge _{g \in S(\m{M})} \,\ell ( ( 1_X - g ) | h | ) \\
      & = \lim _{g \in S(\m{M})} \ell ( ( 1_X - g ) | h | ) \\
      & = \ell ( | h | ) - \lim _{g \in S(\m{M})} \ell ( g | h | ) \\
      & = \ell ( | h | ) - \bigvee _{g \in S(\m{M})} \,\ell ( g | h | ).
\end{align*}
It follows that $q_{\m{M}} (h) = 0$ if and only if
\[  \ell ( | h | ) = \bigvee _{g \in S(\m{M})} \,\ell ( g | h | ), \]
whence the statement.
\end{proof}

\section{Equivalence of $\rho$\,-\,continuity and $\rho$\,-\,regularity}

\begin{proposition}\label{1X}
The functional $\ell$ is $\rho$\,-\,continuous on $\m{M}$
if and only if $S(\m{M})$ contains $1_X$. 
\end{proposition}

\begin{corollary}\label{contreg}
If $\ell$ is $\rho$\,-\,continuous on $\m{M}$,
then $\ell$ is $\rho$\,-\,regular on $\m{M}$.
\end{corollary}

\begin{proof} If $\ell$ is $\rho$\,-\,continuous on $\m{M}$,
then $S(\m{M})$ contains $1_X$ by the preceding
proposition. Then $q_{\m{M}}$ vanishes identically on
$\m{M}$, from which $\ell$ is $\rho$\,-\,regular on $\m{M}$
by virtue of theorem \ref{regular}.
\end{proof}

\begin{proposition}\label{wedge}
For each $I \in \rho$ with $I \subset \m{M}$, one has
\[ \bigwedge _{f \in I} \,\ell (f) = \bigwedge _{f \in I} \,q_{\m{M}} (f). \] 
\end{proposition}

\begin{proof} Let $g \in S(\m{M})$ be arbitrary.
Since $I$ is downward directed, one finds
\begin{align*}
\bigwedge _{f \in I} \,\ell(f) & = \lim _{f \in I} \ell (f)  -  \lim _{f \in I} \ell (gf) \\
                                      & = \lim _{f \in I} \ell ( ( 1_X - g ) f ) \\
                                      & = \bigwedge _{f \in I} \,\ell ( ( 1_X  -  g ) f ). 
\end{align*}
Since $g \in S(\m{M})$ is arbitrary, one also has 
\begin{align*}
\bigwedge _{f \in I} \,\ell (f) & = 
        \bigwedge _{g \in S(\m{M})} \bigwedge _{f \in I} \ell ( ( 1_X - g ) f ) \\
        & = \bigwedge _{f \in I} \bigwedge _{g \in S(\m{M})} \ell ( ( 1_X - g ) f ) \\ 
        & = \bigwedge _{f \in I} \,q_{\m{M}}(f). \qedhere
\end{align*}
\end{proof}

This allows us to reformulate definition \ref{rhocont} in the following way.

\begin{theorem}\label{reform}
The functional $\ell$ is $\rho$\,-\,continuous on $\m{M}$ if and only if
\[ \bigwedge _{f \in I} \,q_{\m{M}}(f) = 0 \textit{ for all }
I \in \rho \textit{ with } I \subset \m{M}. \]
\end{theorem}

\begin{corollary}\label{regcont}
If $\ell$ is $\rho$\,-\,regular on $\m{M}$, then
$\ell$ is $\rho$\,-\,continuous on $\m{M}$. 
\end{corollary}

\begin{proof} If $\ell$ is $\rho$\,-\,regular on $\m{M}$, then
$q_{\m{M}}$ vanishes identically on $\m{M}$ by
theorem \ref{regular}. Theorem \ref{reform}
implies that $\ell$ is $\rho$\,-\,continuous on $\m{M}$.
\end{proof}

\begin{theorem}\label{equiv1}
The functional $\ell$ is $\rho$\,-\,continuous on $\m{M}$
if and only if it is $\rho$\,-\,regular on $\m{M}$.
\end{theorem}

In the light of theorems \ref{equiv1} and \ref{regular}, we can now see
that theorem \ref{reform} is a vast improvement on definition \ref{rhocont}.

\section{The main result}

\begin{definition}\label{K}
We denote the ``kernel'' of the lattice seminorm $q_{\m{M}}$ by
\[ \m{K}(\m{M}) := \{ \,h \in \mathscr{L} : q_{\m{M}}(h) = 0 \,\}. \]
This is a solid vector subspace of $\mathscr{L}$, and thus a vector
sublattice of $\mathscr{L}$. Also, if $\m{M}_1 \subset \m{M}_2$
then $\m{K}(\m{M}_2) \subset \m{K}(\m{M}_1)$.
\end{definition}

\begin{proof}
This follows from the statements in definition \ref{q}.
\end{proof}

\begin{theorem}\label{equiv2}
The functional $\ell$ is $\rho$\,-\,regular on $\m{M}$
if and only if $\m{M} \subset \m{K}(\m{M})$.
\end{theorem}

\begin{proof}
Theorem \ref{regular} and definition \ref{K}.
\end{proof}

\begin{definition}\label{P}
Let $\m{P}(\m{M}) := \m{M} \cap \m{K}(\m{M})$.
This is a vector subspace of $\m{M}$.
\end{definition}

\begin{theorem}\label{regonp}
The functional $\ell$ is $\rho$\,-\,regular on $\m{P}(\m{M})$.
\end{theorem}

\begin{proof}
Put $\m{N} := \m{P}(\m{M}) = \m{M} \cap \m{K}(\m{M})$,
and let $f \in \m{N}$. By theorem \ref{regular} we have
to prove that $q_{\m{N}}(f) = 0$. One one hand, one has
$\m{N} \subset \m{M}$, and so $q_{\m{N}} \leq q_{\m{M}}$
by definition \ref{q}. On the other hand, $f \in \m{K}(\m{M})$,
so that $q_{\m{M}}(f) = 0$ by definition \ref{K}. It follows that
$q_{\m{N}}(f) \leq q_{\m{M}}(f) = 0$.
\end{proof}

\begin{theorem}\label{equiv3}
The functional $\ell$ is $\rho$\,-\,regular on $\m{M}$
precisely when $\m{P}(\m{M}) = \m{M}$.
\end{theorem}

\begin{proof} If $\ell$ is $\rho$\,-\,regular on $\m{M}$, then
$\m{M} \subset \m{K}(\m{M})$ by theorem \ref{equiv2}.
It follows that $\m{P}(\m{M}) = \m{M}$ by definition \ref{P}.
Conversely, if $\m{P}(\m{M}) = \m{M}$, then $\ell$ is
$\rho$\,-\,regular on $\m{M}$ by the preceding theorem
\ref{regonp}.
\end{proof}

\begin{theorem}\label{main1}
The vector subspace $\m{P}(\m{M})$ has the following properties:
\begin{itemize}
\item[$(i)$]  $\ell$ is $\rho$\,-\,continuous on $\m{P}(\m{M})$,
\item[$(ii)$] $\ell$ is $\rho$\,-\,continuous on $\m{M}$
                     if and only if $\m{P}(\m{M})$ is all of $\m{M}$.
\end{itemize}
\end{theorem}

\begin{proof}
Theorems \ref{equiv1}, \ref{regonp}, and \ref{equiv3}.
\end{proof}

The preceding theorem \ref{main1} is our main result. It suggests
that the subspace $\m{P}(\m{M})$ of $\m{M}$ is a ``large'' subspace
of $\rho$\,-\,continuity.

\smallskip
(A largest subspace of $\rho$\,-\,continuity need not exist in the
present generality, as is shown by an argument communicated
to me by Torben Maack Bisgaard.)

\section{Properties of the map $\m{M} \mapsto \m{P}(\m{M})$}

\begin{theorem}
One has $\m{P}(\m{P}(\m{M})) = \m{P}(\m{M})$.
\end{theorem}

\begin{proof}
This follows from theorem \ref{main1} by replacing $\m{M}$ with
$\m{P}(\m{M})$.
\end{proof}

We shall denote by $\| \cdot \|$ the usual lattice seminorm
on $\mathscr{L}$ given by $\| \,f \,\| = \ell ( | f | )$ for all $f \in \mathscr{L}$.

\begin{proposition}\label{closed}
The seminorm $q_{\m{M}}$ is dominated by $\| \cdot \|$, and
thereby is continuous on $( \mathscr{L}, \| \cdot \| )$. It follows that
$\m{K}(\m{M})$ is a closed subspace of $( \mathscr{L}, \| \cdot \| )$.
\end{proposition}

\begin{corollary}\label{csvsl}
The set $\m{K}(\m{M})$ is a closed solid
vector sublattice of $( \mathscr{L}, \| \cdot \| )$.
\end{corollary}

\begin{theorem}\label{inheritance}
The following inheritance properties hold. \\
\noindent If $\m{M}$ is a vector sublattice of $\mathscr{L}$, so is $\m{P}(\m{M})$. \\
\noindent If $\m{M}$ furthermore is Stonean, so is $\m{P}(\m{M})$. \\
\noindent If $\m{M}$ is solid in $\mathscr{L}$, so is $\m{P}(\m{M})$. \\
\noindent If $\m{M}$ is closed in $( \mathscr{L}, \| \cdot \| )$, so is $\m{P}(\m{M})$.
\end{theorem}

\begin{proof}
Definition \ref{P} and the preceding corollary \ref{csvsl}.
\end{proof}

\begin{corollary}\label{properties}
The vector space $\m{P}(\mathscr{L}) = \m{K}(\mathscr{L})$ is the ``kernel'' of the
lattice seminorm $q_{\mathscr{L}}$. It is a closed solid vector sublattice of
$( \mathscr{L}, \| \cdot \| )$. If $\mathscr{L}$ is Stonean, so is $\m{P}(\mathscr{L})$.
\end{corollary}

\begin{proof} Definitions \ref{P} and \ref{K},
corollary \ref{csvsl} and theorem \ref{inheritance}.
\end{proof}

\begin{theorem}
If $\m{M}_1 \subset \m{M}_2$, then
$\m{P}(\m{M}_2) \cap \m{M}_1 \subset \m{P}(\m{M}_1)$.
\end{theorem}

\begin{proof} One has
$\m{P}(\m{M}_2) \cap \m{M}_1 \subset \m{K}(\m{M}_2) \cap \m{M}_1
\subset \m{K}(\m{M}_1) \cap \m{M}_1 = \m{P}(\m{M}_1)$ by
definition \ref{P} and the last statement in definition \ref{K}.
\end{proof}

Whence, as special cases, the following three corollaries:

\begin{corollary}
One has $\m{P}(\mathscr{L}) \cap \m{M} \subset \m{P}(\m{M})$.
\end{corollary}

\begin{corollary}
If $\m{M}_1$ is a subspace between $\m{P}(\m{M}_2)$ and
$\m{M}_2$, that is, if $\m{P}(\m{M}_2) \subset \m{M}_1 \subset \m{M}_2$,
then $\m{P}(\m{M}_2) \subset \m{P}(\m{M}_1) \subset \m{M}_1 \subset \m{M}_2$.
\end{corollary}

\begin{corollary}
If $\m{P}(\mathscr{L}) \subset \m{M}$, then
$\m{P}(\mathscr{L}) \subset \m{P}(\m{M}) \subset \m{M}$.
\end{corollary}

\end{document}